\documentclass{article}

\usepackage[utf8]{inputenc}
\usepackage[T1]{fontenc}
\usepackage{geometry}
\usepackage{amsmath, amsthm, amssymb}

\newtheorem{thm}{Theorem}[section]
\newtheorem{cor}[thm]{Corollary}
\newtheorem{lem}[thm]{Lemma}

\newtheorem{rmk}[thm]{Remark}
\newtheorem{definition}[thm]{Definition}

\title{Local neighborliness of the symmetric moment curve}
\author{Seung Jin Lee}
\date{}

\begin{document}
\maketitle
\centerline {\bf Abstract}
A centrally symmetric analogue of the cyclic polytope, the bicyclic polytope, was defined in [BN08]. The bicyclic polytope is defined by the convex hull of finitely many points on the symmetric moment curve where the set of points has a symmetry about the origin. In this paper, we study the Barvinok-Novik orbitope, the convex hull of the symmetric moment curve.
 It was proven in [BN08] that the orbitope is locally $k$-neighborly, that is, the convex hull of any set of $k$ distinct points on an arc of length not exceeding $\phi_k$ in $\mathbb{S}^1$ is a $(k-1)$-dimensional face of the orbitope for some positive constant $\phi_k$. We prove that we can choose $\phi_k $ bigger than $\gamma k^{-3/2} $ for some positive constant $\gamma$.

\section{Introduction and main result}

Let $P$ be a $d$-dimensional polytope with $n$ vertices and let $f_i(P)$ be the number of $i$-dimensional faces of $P$. It is known that $f_j(P) \leq  \left(\begin{smallmatrix} n\\ j+1 \end{smallmatrix}\right)$ for $j \leq \lfloor \frac{d}{2}\rfloor$ and equality holds when $P$ is the cyclic polytope. However, the situation for centrally symmetric polytopes is different. For example, the largest number of edges, fmax$(d,n;1)$, that a $d$-dimensional centrally symmetric polytope on $n$ vertices can have is unknown even for $d=4$. In [BN08], for fixed even dimension $d=2k$ and an integer $1\leq j <k$ Barvinok and Novik proved that $\textrm{fmax}(d,n;j)$, the maximum number of $j$-dimensional faces of a centrally symmetric $d$-dimensional polytope with $n$ vertices, is at least $(c_j(d)+o(1)) \left(\begin{smallmatrix} n\\ j+1 \end{smallmatrix}\right)$ for some $c_j(d)>0$ and at most $(1-2^{-d}+o(1)) \left(\begin{smallmatrix} n\\ j+1 \end{smallmatrix}\right)$ as $n$ grows. The authors also proved that $c_1(d)\geq 1-\frac{1}{d-1}$ and $c_j(d)>0$ for any $j\leq k-1$. To get a lower bound we need to define a centrally symmetric analog of cyclic polytopes - bicyclic polytopes.
\bigskip

As in [BN08], the authors consider the convex hull of the \emph{symmetric moment curve} \\
\[ SM_{2k}(t)=(\cos t, \sin t, \cos 3t, \sin 3t, \ldots, \cos(2k-1)t, \sin(2k-1)t)\]
This curve is centrally symmetric. We define the Barvinok-Novik orbitope
\[{\cal{B}}_{2k}  = \textrm{conv}(SM_{2k}(t): 0 \leq t \leq 2\pi )\]

In [BN08], it is proven that  ${\cal B}_{2k}$ is a locally $k$-neighborly.

\begin{thm}
For every positive integer $ k $ there exists a number $\psi_k > 0$ such that if $t_1,\ldots,t_k \in S^1$ are distinct points that lie on an arc of length less than $\psi_k$, then \\
\centerline{\textrm{conv}$(SM_{2k}(t_1), \ldots , SM_{2k}(t_{k}))$}
is a (k-1)-dimensional face of ${\cal B} _{2k}$.
\end{thm}

In this paper, a \emph{face} of a convex body is an exposed face, intersection of the body with a supporting hyperplane. Let $\phi_k$ be the supremum of all possible value of $\psi_k$ in Theorem 1.1. Then $\phi_k$ also satisfies Theorem 1.1 because if $t_i$'s are distinct points in $\mathbb{S}^1$ lying on an arc of length less than $\phi_k$, the points also lie on an arc of length $\psi $ less than $\phi_k$ for some $\psi$ satisfying Theorem 1.1.

\bigskip

The goal of this paper is to find a lower bound of $\phi_k$. 

\begin{thm}
Let $\phi_k$ be the supremum of all possible value of $\psi_k$ in Theorem 1.1. Then we have $\phi_k>\sqrt{6}k^{-3/2}$.
\end{thm}

The idea of the proof of Theorem 1.2 is as follows. First of all, we find a lower bound of a distance between the boundary of the Barvinok-Novik orbitope and the origin by studying the minimum volume ellipsoid of the Barvinok-Novik orbitope. Then we show that for any $k$ points lying on an arc of length less than $\sqrt{6}k^{-3/2}$, there is no new intersection point between the affine hyperplane that is tangent to the symmetric moment curve at the points and the Barvinok-Novik orbitope. To be more precise, if arc of length is too small then there is no new point on the opposite arc because the line segment joining new point and one of $k$ points will pass through the interior of the Barvinok-Novik orbitope. However, such new point may appear only on the opposite arc , as shown in [BN08], because the hyperplane is a supporting hyperplane of ${\cal B}_{2k}$.\\

In Section 2 we discuss the minimum volume ellipsoid of the Barvinok-Novik orbitope and find a lower bound of a distance between the origin and the boundary of the Barvinok-Novik orbitope. In Section 3 we prove Theorem 1.2. 

\section{The minimum volume ellipsoid of the Barvinok-Novik orbitope}

\smallskip

In this section, we prove that ${\cal{B}}_{2k}$ contains the sphere of radius $\frac{1}{\sqrt{2}}$ centered at the origin. To prove this result, we need the notion of the minimum volume ellipsoid and two theorems related to it.

\begin{definition}
Given a convex body (compact convex set with a non-empty interior) $B \subset \mathbb{R}^d$, there is a unique ellipsoid $E_{\min} \supset B$ of the minimum volume, called the minimum volume ellipsoid of $B$. (See, for example, $[\textrm{B}97]$)

\end{definition}
\begin{thm} $[\textrm{BB} 05]$
Let $G$ be a compact group acting on the Euclidean space $V$ with $G$-invariant inner product $\langle\cdot,\cdot\rangle$ and $v$ be a nonzero vector in $V$. Let $B$ be the convex hull of the orbit of a vector $v\in V$:
\[ B= \textrm{conv}(gv:g\in G).\]
Suppose that the affine hull of B is V .
Then there exists a decomposition 
\[V=\bigoplus _i V_i\]
of V into the direct sum of pairwise orthogonal irreducible components such that
the following holds.\\

The minimum volume ellipsoid $E_{min}$ of $B$ is defined by the inequality
\[E_{\min} =\left\{x :\quad \sum_i \frac{\dim V_i}{\dim V}\cdot \frac{\langle x_i,x_i\rangle }{\langle v_i,v_i \rangle} \leq 1\right\},\]
where $x_i$ (resp. $v_i$) is the orthogonal projection of $x$ (resp. $v$) onto $V_i$.
\end{thm}
\begin{thm} If a convex body $B$ is symmetric about the origin, then $(\dim B)^{-1/2} E_{\min} \subset B \subset E_{\min}$. (See, for example, $[\textrm{B}97]$)
\end{thm}
\bigskip

In our situation, we have a following corollary.

\begin{cor}

$B$ contains the sphere of radius $\frac{1}{\sqrt{2}}$.
\end{cor}

\noindent{\bf Proof of Corollary 2.4}. 
In our case, we have $V=\mathbb R^{2k}$, $v=(1,0,1,0,\ldots,1,0)$, $G=\mathbb{S}^1$ and the action of $G$	 is given by the matrix
\[ \left(\begin{matrix} \cos(t) & -\sin(t) & 0& 0& \cdots & 0 & 0 
\\ \sin(t)&\cos(t) & 0& 0& \cdots & 0 & 0 
\\ 0 & 0 & \cos(3t) &-\sin(3t) & \cdots & 0 & 0
\\ 0& 0 & \sin(3t) & \cos(3t) & \cdots & 0 & 0
\\ \vdots & \vdots & \vdots &\vdots &\ddots & \vdots & \vdots
\\ 0& 0 & 0& 0 &\cdots & \cos((2k-1)t) & -\sin((2k-1)t)
\\ 0& 0 & 0& 0 &\cdots & \sin((2k-1)t) & \cos((2k-1)t)
\end{matrix}\right).\] 
In particular, the decomposition of $V=\bigoplus_i V_i$ is multiplicity-free because $V_i$' s are $\mathbb R^2$ with the action of $G$ by multiplication of 
 \[ \left(\begin{matrix} \cos((2j-1)t) & -\sin((2j-1)t)\\ \sin((2j-1)t)&\cos((2j-1)t) \end{matrix}\right)\] 
so $V_i$' s are not isomorphic each other. Since $v=(1,0,1,0,\ldots,1,0)$ we have $\langle v_i,v_i \rangle=1$ for any $i$. Therefore, $E_{\min}$ is the sphere of radius $\sqrt{k}$ centered at the origin by Theorem 2.2. By Theorem 2.3, $B$ contains the sphere of radius $\frac{1}{\sqrt{2}}$. $\qed$

\section{Proof of Theorem 1.2.}

In this section, we prove that $\phi_k$ is bigger than $\sqrt{6} k^{-3/2}$. For simplicity, we use $x(t)$ instead of $SM_{2k}(t)$.
\smallskip

\begin{thm} 
Suppose that there exist a positive number $\psi$ less than $\pi$, positive even integers $m_i$ for $i=1,\ldots,l$ satisfying $\sum_{i=1}^{l}m_i=2k$ and $l$ points $t_1,\ldots,t_l$ on an arc of length $\psi$ such that the affine hyperplane $H$ tangent to $x(t)$ at each point $t_i$ with multiplicity $m_i$ is a supporting hyperplane of ${\cal B}_{2k}$ and intersects with an opposite arc at another point $x(s)$. Then
\[ |s-\pi-t_i| > \sqrt{3/2}k^{-3/2} \textrm{ for any } i.\]
\end{thm}
\smallskip

{\bf Proof of Theorem 3.1}. 
Suppose that $-\frac{\psi}{2} \leq t_i \leq \frac{\psi}{2}$ for any $i$. It is known that new point $x(s)$ should lie in the opposite arc (See [BN08] Lemma 6.3), so we can assume that $\pi-\frac{\psi}{2} \leq s \leq \pi + \frac{\psi}{2}$. In this situation, we prove that distance between $x(s)$ and an opposite point of any $x(t_i)$ cannot be too small.

Since $x(s)$ and $x(t_i)$ are vertices of the face defined by $x(t_i)$'s, the midpoint $\frac{x(s)+x(t_i)}{2}$ lies on the face. Therefore, we have $\big|\frac{x(s)+x(t_i)}{2}\big|^2 \geq \frac{1}{2} $ because ${\cal{B}}_{2k}$ contains a sphere of radius $\frac{1}{\sqrt{2}}$ centered at the origin. \\

If $|s-\pi-t_i|=\epsilon$, we have
\[  \Big|  \frac{x(s)+x(t_i)}{2}\Big|^2-\frac{1}{2}=-\frac{1}{2}+\frac{1}{4}\sum_{i=1}^k\left(\cos((2i-1)s)+\cos((2i-1)t_i))^2+(\sin((2i-1)s)+\sin((2i-1)t_i)\right)^2\]
\[=-\frac{1}{2} +\frac{1}{2}\sum_{i=1}^k\left( 1+ \cos((2i-1)(s-t_i))\right )=\frac{1}{2}\left( k-1 + \frac{\sin(2k(s-t_i))}{2\sin(s-t_i)} \right) \]
\[=\frac{1}{2}\left( k-1 - \frac{\sin(2k\epsilon)}{2\sin\epsilon}\right )\]

Now we are using following well-known inequalities.
\[x-x^3/6<\sin(x)<x  \textrm{ \hspace{0.5cm}     for  } x>0\]

Therefore, we have

\[\frac{1}{2}\left( k-1 - \frac{\sin(2k\epsilon)}{2\sin\epsilon}\right )< \frac{1}{2}\left( k-1- \frac{2k\epsilon-8k^3\epsilon^3/6}{2\epsilon}\right )= -1/2 + k^3\epsilon^2/3\] for $\epsilon>0$.
Since this is greater than or equal to zero, we have $\epsilon>\sqrt{3/2}k^{-3/2}$. \qed

\smallskip

We need one more lemma to prove Theorem 1.2. \\
\smallskip
\begin{lem}
Let $\Gamma \subset \mathbb{S}^1$ be an arc of length less than $\pi$. Let $t_i\subset \Gamma, i=1,\ldots,l$ be distinct points and let $m_i>1, i=1,\ldots,l$ be integers such that
\[\sum_{i=1}^l m_i=2k\]
Then the following $2k$ vectors
\[x(t_i)-x(t_l) \textrm{ for } i=1,\ldots, l-1,\]
\[\frac{d^n}{dt^n} x(t)\Big|_{t=t_i} \textrm{ for } n=1,\ldots,m_i -1 \textrm{ and } i=1,\ldots,l \]
\[\frac{d^{m_1}}{dt^{m_1}} x(t)\Big|_{t=t_1}\]
are linearly independent in $\mathbb{R}^{2k}$

In particular, there exists a unique affine hyperplane $H\subset \mathbb{R}^{2k}$ that is tangent to $x(t)$ at each point $t_i$ with multiplicity $m_i$.
\end{lem}

{\bf Proof.} Assume that the vectors are not linearly independent. Then there exists a non-zero vector $a\in \mathbb{R}^{2k}$ which is orthogonal to all the vectors. Let us define a trigonometric polynomial
\[p(t)=\langle a, x(t)-x(t_l)\rangle,\]
Hence $p(t)$ is not identically zero and has zeroes at $t_i$ with multiplicity $m_i$ respectively for $i=2,\ldots, l$ and a zero at $t_1$ with multiplicity $m_1+1$. Therefore the total number of roots of $p(t)$ on $\Gamma$, counting multiplicities, is at least $2k+1$. By Rolle's Theorem, the number of roots of the derivative $p'(t)$ on $\Gamma$ is at least $2k$, counting multiplicities. However, the constant term of $p'(t)$ is 0, so we have $p'(t+\pi)=-p'(t)$ and the total number of roots of $p'$ on the circle is at least $4k$, counting multiplicities. However, since $p'(t)$ is a trigonometric polynomial of degree $2k-1$ it has at most $4k-2$ roots if $p'(t)$ is nonzero. Hence $p'(t)\equiv 0$ and $p(t)$ is a constant, which is a contradiction.
\qed
\\
\smallskip

{\bf Proof of Theorem 1.2.}  Let $t_1\leq t_2\leq\ldots\leq t_k$ be $k$ points on an arc $\Gamma$ of length $\phi_k$ of $\mathbb{S}^1$.\\

Let us define a function
\[p(t_1,\ldots,t_k)=\textrm{dist}(H_{t_1,t_2,\ldots,t_k},x(\Gamma+\pi))\]
where $H_{t_1,t_2,\ldots,t_k}$ is the affine hyperplane that is tangent to $x(t)$ at each point $t_i$ with even multiplicity $m_i$. Note that we use multiset notation such that number of $t_i$ appeared in the multiset is same as $m_i/2$, so we know that sum of all the $m_i$'s are $2k$. Note that $H_{t_1,t_2,\ldots,t_k}$ is well-defined because of Lemma 3.2. \\

Now we take an infimum of $p(t_1,\ldots,t_k)$ for all $t_i$'s lying on an arc of length at most $\phi_k$ centered at 0. By definition of $\phi_k$, the infinum is nonnegative. If the infimum is strictly positive, by continuity of $p$ we can extend the length of the arc bigger than $\phi_k$ such that the infimum of $p$ for $k$ points lying on the bigger arc is positive. However, this means that for any $k$ distinct points lying on the bigger arc the affine hyperplane tangent at these points is a supporting hyperplane, contradicting definition of $\phi_k$.\\

Hence the infimum is 0. Since domain of the function $p$ is compact, we can find $l$ points $t_1,\ldots,t_l$ on $\Gamma$ and positive even integers $m_i$ for $i=1,\ldots,l$ satisfying $\sum_{i=1}^{l}m_i=2k$ such that the affine hyperplane $H$ is  a supporting hyperplane of ${\cal B}_{2k}$ tangent to $x(t)$ at each point $t_i$ with multiplicity $m_i$ and intersects with the opposite arc $x(\Gamma+\pi)$ at some point, say $x(s)$. If $\phi_k$ is less than or equal to $\sqrt{6}k^{-3/2}$, there exists a point $t_i$ such that $ |s-\pi-t_i| \leq \sqrt{3/2}k^{-3/2}$ and it contradicts with Theorem 3.1. \qed

\section{Concluding remarks}

\begin{rmk}
The estimate of the Corollary 2.2 can be improved by $\Omega(k^{-5/4})$ from $\Omega(k^{-3/2})$. This can be done by considering a linear combination of $-x(s)$ and two points among $x(t_i)$' s close to $-x(s)$ (say, $x(t_i)$ and $x(t_j)$). To be more precise, we can get a better bound by considering
\[ \left|\frac{x(s)}{2}+ \frac{(s-t_j)x(t_i)}{2(t_i-t_j)} + \frac{(t_i-s)x(t_j)}{2(t_i-t_j)} \right|\]

\end{rmk}
\begin{rmk}
The Barvinok-Novik orbitope ${\cal B}_{2k}$ does not contains a sphere of radius bigger than 1. In fact, a hyperplane $x_{2k-1}=1$ defines $(2k-2)$-dimensional face of ${\cal B}_{2k}$ and the distance between the hyperplane and the origin is 1.
\end{rmk}

{\bf Acknowledgements}
Thanks to Alexander Barvinok for many helpful discussions. This research was partially supported by NSF grant DMS 0856640.

\bigskip

Department of Mathematics, University of Michigan, Ann Arbor, MI 48109-1043, USA\\
\emph{E-mail address}: lsjin@umich.edu

\end{document}